\documentclass[12pt]{amsart}

\usepackage[latin1]{inputenc}
\usepackage{amsmath}
\usepackage{amsfonts}
\usepackage{amssymb}
\usepackage{graphics}
\usepackage{enumerate}
\usepackage{amssymb,amsmath,amsthm,amscd,latexsym,verbatim,graphicx,amsfonts}

\newtheorem{theorem}{Theorem}

\newtheorem{question}{Question}
\newtheorem{conjecture}{Conjecture}

\newcommand{\bbR}{\mathbb{R}}

\newcommand{\bbN}{\mathbb{N}}

\newcommand{\nri}{n\rightarrow\infty}

\begin{document}
\title[ ] {On Some Convexity Questions of Handelman}

\bibliographystyle{plain}

\thanks{  }


\maketitle

\begin{center}
\textbf{Brian Simanek}
\end{center}


\begin{abstract}
We resolve some questions posed by Handelman in 1996 concerning log convex $L^1$ functions.  In particular, we give a negative answer to a question he posed concerning the integrability of $h^2(x)/h(2x)$ when $h$ is $L^1$ and log convex and $h(n)^{1/n}\rightarrow1$.
\end{abstract}

\vspace{5mm}

\section{Introduction}\label{intro}

In \cite{Hand}, Handelman investigated eventual positivity of power series and deduced its existence for a wide variety of functions by appealing to a particular maximal function.  If $h:(0,\infty)\rightarrow\infty)$ is continuous, then he defined the maximal function
\[
H_h(a)=\max_{b\geq a}\frac{h(b)}{h(a+b)}.
\]
This maximal function was introduced in \cite{Hand}, where some of its properties are discussed.  In particular, it is meaningful to have an understanding of when $hH_h$ is integrable on $(0,\infty)$.  It is easy to see that if $h$ is log convex, then $H_h(x)=h(x)/h(2x)$.  This led Handelman to ask  the following question (see \cite[page 338]{Hand}):

\begin{question}\label{q1}  If $h:(0,\infty)\rightarrow(0,\infty)$ is a log convex function that is integrable on $(0,\infty)$ and satisfies $\lim_{\nri}h(n)^{1/n}=1$, then is it true that $h(x)^2/h(2x)$ is also integrable on $(0,\infty)$?
\end{question}

One of our main results is a demonstration that the answer to Question \ref{q1} is ``no."  In fact, we will prove the following result:

\begin{theorem}\label{allr}
There is a function $h:(0,\infty)\rightarrow(0,\infty)$ that is log convex, integrable on $(0,\infty)$, satisfies $\lim_{\nri}h(n)^{1/n}=1$, and is such that $h^{r}(x)/h(rx)$ is not integrable on $(0,\infty)$ for any $r>0$.
\end{theorem}

\noindent Our proof of Theorem \ref{allr} is constructive in that we will show how to actually create a counterexample.

Also in \cite{Hand}, Handelman made the following conjecture (see the discussion following \cite[Theorem 9]{Hand}):

\begin{conjecture}\label{conj1}
Suppose $h:(0,\infty)\rightarrow(0,\infty)$ is a log convex function that satisfies
\begin{align}\label{hrat}
\lim_{\nri}\frac{h(n)}{h(n+1)}=1
\end{align}
and $h^2(x)/h(2x)$ is integrable on $(0,\infty)$.  Then $h$ is also integrable on $(0,\infty)$.
\end{conjecture}

Our second main result is the following:

\begin{theorem}\label{t1}
Conjecture \ref{conj1} is true.
\end{theorem}

\noindent\textit{Remark.}  We should point out that the hypothesis (\ref{hrat}) is essential to proving Theorem \ref{t1}, for otherwise one could take $h(x)=x^x$ as a counterexample.

\medskip

The remainder of the paper is devoted to the proofs of Theorem \ref{allr} and Theorem \ref{t1}.  Our methods are elementary and require only basic facts about convex functions (see \cite{SimConv} for a discussion of many tools in convexity theory).

\section{Proofs}\label{proof}

The purpose of this section is to prove all of the results discussed in the introduction.

\subsection{The Construction}\label{constr}

In this section, we will resolve Question \ref{q1}.  Let us write $h(x)=\exp(f(x))$ where $f(x)$ is convex.  Since $h\in L^1(\bbR^+)$, it must be that $\lim_{x\rightarrow\infty}f(x)=-\infty$.  The remaining condition on $h$ implies $\lim_{\nri}f(n)/n=0$, or equivalently (by the convexity of $f$) $\lim_{x\rightarrow\infty}f'(x)=0$ provided $f'(x)$ exists.  In fact, the function $f$ we construct will be piecewise linear and hence $f'(x)$ will be undefined on a discrete set.  We will choose sequences $\{a_n\}_{n=0}^{\infty}$, $\{m_n\}_{n=0}^{\infty}$, and $\{b_n\}_{n=0}^{\infty}$ so that
\begin{equation}\label{fdef}
f(x)=m_nx+b_n\qquad,\qquad x\in[a_n,a_{n+1}],
\end{equation}
and $f$ is continuous.

To begin our construction, let $\{m_n\}_{n=0}^{\infty}$ be a fixed sequence of negative real numbers that monotonically increases to $0$.  With this fixed sequence in hand, we will construct the sequence $\{a_n\}_{n=0}^{\infty}$ inductively, and the sequence $\{b_n\}_{n=0}^{\infty}$ will then be defined implicitly in order to make $f$ continuous.

We begin our construction of the sequence $\{a_n\}_{n\geq0}$ by defining $a_0=0$ and we also define $b_0=0$.  Now, choose $a_1$ large enough so that
\begin{itemize}
\item
\[
-\frac{e^{m_0a_1}}{m_1}<\frac{1}{2}
\]
\item
\[
a_1>1=a_0+e^{-b_0}.
\]
\end{itemize}
Now set $b_1=m_0a_1-m_1a_1$ and observe that $m_1a_1+b_1=m_0a_1$.

Now let us assume that $\{a_j\}_{j=0}^n$ and $\{b_j\}_{j=0}^n$ have already been defined.  We will now show how one can choose $a_{n+1}$ and then $b_{n+1}$ to complete the construction.  Indeed, as above, we will choose $a_{n+1}$ large enough so that
\begin{align*}
-\frac{e^{m_na_{n+1}+b_n}}{m_{n+1}}<\frac{1}{2^{n+1}},
\end{align*}
\begin{align*}
a_{n+1}>\max\left\{a_n+e^{(1-t)b_n}:t\in(0,n]\right\}.
\end{align*}
Then define
\[
b_{n+1}=m_na_{n+1}+b_n-m_{n+1}a_{n+1}
\]
and observe that $m_{n+1}a_{n+1}+b_{n+1}=m_na_{n+1}+b_n$. Proceeding inductively, we arrive at two sequences $\{a_n\}_{n=0}^{\infty}$ and $\{b_n\}_{n=0}^{\infty}$.  It is clear from our construction that $a_{n+1}>a_n+1$ (since $1\in(0,n]$) and so $\lim_{\nri}a_n=\infty$.  Therefore this procedure defines $f$ on all of $(0,\infty)$ if we define $f$ by \eqref{fdef}.

Now let us check that this function has the desired properties.  First of all, since $m_n\rightarrow0$ monotonically, it is clear that $f(n)/n\rightarrow0$ and also that $f$ is convex.  Now we calculate
\begin{align*}
\int_0^{\infty}h(x)dx&=\int_0^{\infty}e^{f(x)}dx=\sum_{n=0}^{\infty}\int_{a_n}^{a_{n+1}}e^{m_nx+b_n}dx\\
&\leq\sum_{n=0}^{\infty}\int_{a_n}^{\infty}e^{m_nx+b_n}dx\\
&=\sum_{n=0}^{\infty}-\frac{e^{m_na_n+b_n}}{m_n}\\
&=-\frac{1}{m_0}+\sum_{n=1}^{\infty}-\frac{e^{m_{n-1}a_n+b_{n-1}}}{m_n}\qquad\qquad\mbox{(we use continuity of $f$ here)}
\\&<-\frac{1}{m_0}+\sum_{n=1}^{\infty}2^{-n},
\end{align*}
which is clearly finite.
Therefore, $h\in L^1(\bbR^+)$ as desired.

Finally, fix $r\in(0,\infty)$ and choose $N\in\bbN$ so that $r<N$.  Notice that
\[
\frac{h^r(x)}{h(rx)}=e^{(r-1)b_n}\qquad,\qquad x\in[a_n,a_{n+1}].
\]
Therefore,
\begin{align*}
\int_0^{\infty}\frac{h^r(x)}{h(rx)}dx&=\sum_{n=0}^{\infty}e^{(r-1)b_n}(a_{n+1}-a_n)>\sum_{n=N}^{\infty}e^{(r-1)b_n}(a_{n+1}-a_n)>\sum_{n=N}^{\infty}1,
\end{align*}
by construction, so $\frac{h^r(x)}{h(rx)}\not\in L^1(\bbR^+)$.  This completely answers Question \ref{q1}.

\subsection{The Conjecture}\label{conjans}

In this section, we will prove Theorem \ref{t1}.  The log convexity of $h$ implies that $h$ is either monotone increasing on $(A,\infty)$ for some $A\geq0$ or monotone decreasing on $(0,\infty)$.  In the latter case, we have $h(x)\geq h(2x)$ and so
\[
\int_0^{\infty}h(x)dx\leq\int_0^{\infty}\frac{h^2(x)}{h(2x)}dx<\infty,
\]
so $h$ is integrable.

If $h$ is monotone increasing on $(A,\infty)$, then $\log(h(x))$ is also increasing on $(A,\infty)$.  Since $\log(h(x))$ is convex, it must be that there is some constant $c>0$ so that
\begin{align*}
\log(h(n+1))-\log(h(n))\geq c,\qquad n>A.
\end{align*}
This implies $h(n+1)/h(n)\geq e^c$, which means $h$ cannot satisfy (\ref{hrat}).  Therefore, this case cannot occur, and we have proven Theorem \ref{t1}.

\bigskip

\noindent\textbf{Acknowledgments.}  The author gratefully acknowledges support from the Simons Foundation through collaboration grant 707882.


\vspace{10mm}



\begin{thebibliography}{14}

\bibitem{Hand} D. Handelman, {\em Eventual positivity for analytic functions}, Math. Ann., 304 (1996), no. 2, 315--338.

\bibitem{SimConv} B. Simon, {\em Convexity: An Analytic Viewpoint}, Cambridge Tracts in Mathematics, 187. Cambridge University Press 2011.



\end{thebibliography}
\end{document}